\input amstex\documentstyle{amsppt}  
\pagewidth{12.5cm}\pageheight{19cm}\magnification\magstep1  
\topmatter
\title Partial flag manifolds over a semifield\endtitle
\author G. Lusztig\endauthor
\address{Department of Mathematics, M.I.T., Cambridge, MA 02139}\endaddress
\thanks{Supported by NSF grant DMS-1855773.}\endthanks
\endtopmatter   
\document

\define\mpb{\medpagebreak}

\define\sqc{\sqcup}

\define\qua{\quad}

\define\op{\oplus}
   
\define\part{\partial}
\define\emp{\emptyset}

\define\m{\mapsto}
\define\do{\dots}

\define\sub{\subset}    

\define\T{\times}

\define\nl{\newline}
\redefine\i{^{-1}}

\define\un{\underline}

\define\ot{\otimes}

\define\Hom{\text{\rm Hom}}
\define\End{\text{\rm End}}

\define\supp{\text{\rm supp}}

\redefine\b{\beta}

\redefine\o{\omega}

\redefine\t{\tau}
\define\th{\theta}

\redefine\l{\lambda}

\define\x{\xi}

\redefine\G{\Gamma}

\define\CC{\bold C}

\define\NN{\bold N}

\define\RR{\bold R}

\define\ZZ{\bold Z}

\define\cc{\Cal C}
\define\cd{\Cal D}

\define\cp{\Cal P}

\define\cs{\Cal S}

\define\cv{\Cal V}

\define\cx{\Cal X}

\define\fG{\frak G}

\define\bJ{\bar J}
\define\bE{\bar E}

\head Introduction\endhead
\subhead 0.1\endsubhead
Let $G$ be the group of simply connected type associated in \cite{MT}, \cite{Ma}, \cite{Ti}, \cite{PK},
 to a not necessarily positive definite symmetric Cartan matrix and to the field $\CC$. We 
assume that a pinning of $G$ is given. It consists of a ``Borel subgroup'' $B^+$, a ``maximal torus''
$T\sub B^+$ and one parameter subgroups $x_i:\CC@>>>G,y_i:\CC@>>>G\qua (i\in I)$ analogous to those
in \cite{L94}. We have $x_i(\CC)\sub B^+$.
We fix a subset $J\sub I$. Let $\Pi^J$ be the subgroup of $G$ generated by $B^+$ and by 
$\cup_{i\in J}y_i(\CC)$. Let $\cp^J$ be the set of subgroups of $G$ which are $G$-conjugate to
$\Pi^J$ (a partial flag manifold). As in \cite{L94, 2.20} we consider the submonoid $G_{\ge0}$ of $G$
generated by $x_i(a),y_i(a)$ with $i\in I,a\in\RR_{\ge0}$ and by 
the ``vector part'' $T_{>0}$ of $T$. ($T$ is a product of $T_{>0}$ and a compact torus.)
Let $K$ be a semifield. 
Let $\fG(K)$ be the monoid associated to $G,K$ by generators and relations in \cite{L18, 3.1(i)-(viii)}.
When $K=\RR_{>0}$ this can be identified with $G_{\ge0}$ by an argument given in \cite{L19a}.
In this paper we define an analogue of $\cp^J$ in the case where $\CC$ is replaced by any semifield $K$. 
This is a set $\cp^J(K)$ with an action of $\fG(K)$. 

A part of our argument involves a construction of an analogue of the highest weight integrable
representations of $G$ when $G$ is replaced by the monoid $\fG(K)$. The possibility of such a
construction comes from the positivity properties of the canonical basis \cite{L93}.

\subhead 0.2\endsubhead
In this subsection we assume that our Cartan matrix is of finite type.
If $K=\RR_{>0}$, the set $\cp^J(K)$ coincides with the subset $\cp^J_{\ge0}$ of $\cp^J$ defined
in \cite{L98}. If $K$ is the semifield $\ZZ$ and $J=\emp$, a definition of the flag manifold over $\ZZ$
was given in \cite{L19b}; we expect that it agrees with the definition in this paper, but we have not
proved that.
In the case where $G=SL_n$, a form over $\ZZ$ of a grassmannian was defined earlier in \cite{SW}.

\head 1. The set $\cp^J(K)$\endhead
\subhead 1.1\endsubhead
Let $\cx=\Hom(T,\CC^*)$. This is a free abelian group with basis $\{\o_i;i\in I\}$ consisting of 
fundamental weights. Let $\cx^+=\sum_{i\in I}\NN\o_i\sub\cx$ be the set of dominant weights.
For $\l\in\cx$ let $\supp(\l)$ be the set of all $i\in I$ such that $\o_i$ appears with $\ne0$ coefficient
in $\l$. Let $\cx^+_J=\{\l\in\cx^+;\supp(\l)=I-J\}$, $\cx^+_{\bJ}=\{\l\in\cx^+;\supp(\l)\sub I-J\}$.

The irreducible highest weight integrable representations of $G$ are indexed by their highest weight, an 
elements of $\cx^+$. For $\l\in\cx^+$ let ${}^\l V$ be a $\CC$-vector space which is an irreducible 
highest weight integrable representation of $G$ indexed by $\l$. 
Let ${}^\l P$ be the set of lines in ${}^\l V$.
Let ${}^\l\x^+$ be a highest weight vector 
of ${}^\l V$. 
Let ${}^\l\b$ be the canonical basis of ${}^\l V$ (see \cite{L93, 11.10}) containing ${}^\l\x^+$. 

For a nonzero vector $\x$ in a vector space $V$ we denote by $[\x]$ the line in $V$ that contains $\x$. 
Note that $\Pi^J$ (see 0.1) is the stabilizer of $[{}^\l\x]$ in $G$ where $\l\in\cx^+_J$. 

For $\l,\l'$ in $\cx^+$ we define a linear map 
$$E:{}^\l V\T{}^{\l'}V@>>>{}^\l V\ot{}^{\l'}V$$ by 
$(\x,\x')\m\x\ot\x'$ and a linear map
$$\G:{}^{\l+\l'}V@>>>{}^\l V\ot{}^{\l'}V$$
which is compatible with the $G$-actions and takes ${}^{\l+\l'}\x^+$ to ${}^\l\x^+\ot{}^{\l'}\x^+$.
Let ${}^{\l,\l'}P$ be the set of lines in ${}^\l\x^+\ot{}^{\l'}\x^+$. Now $E$ induces a map
$\bE:{}^\l P\T{}^{\l'}P@>>>{}^{\l,\l'}P$  and $\G$ induces a map $\bar\G:{}^{\l+\l'}P@>>>{}^{\l,\l'}P$.

Let $\cc$ be the set of all collections $\{x_\l\in{}^\l V;\l\in\cx^+_{\bJ}\}$
such that for any $\l,\l'$ in $\cx^+_{\bJ}$ we have $\G(x_{\l+\l'})=E(x_\l,x_{\l'})$.
Let $\cc^*$ be the set of all $(x_\l)\in\cc$ such that $x_\l\ne0$ for any $\l\in\cx^+_{\bJ}$.
Let $H$ be the group consisting of all collections $\{z_\l\in\CC^*;\l\in\cx^+_{\bJ}\}$
such that for any $\l,\l'$ in $\cx^+_{\bJ}$ we have $z_{\l+\l'}=z_\l z_{\l'}$.
Now $H$ acts on $\cc$ by $(z_\l),(x_\l)\m(z_\l x_\l)$.
This restricts to a free action of $H$ on $\cc^*$. Let ${}'\cp^J$ be the set of orbits for this action.
Note that $G$ acts on $\cc$ by $g(x_\l)=(g(x_\l))$. This induces a $G$-action on $\cc^*$ and on ${}'\cp^J$.
We define a map $\th:\cp^J@>>>{}'\cp^J$ by $g\Pi^Jg\i\m\text{$H$-orbit of }(g({}^\l\x))$ where $g\in G$.
This is well defined since $({}^\l\x)\in\cc$ and since for $g\in\Pi^J$, 
$(g({}^\l\x))$ is in the same $H$-orbit as $({}^\l\x)$. We show:

\proclaim{Lemma 1.2} $\th:\cp^J@>>>{}'\cp^J$ is a bijection.
\endproclaim
For $\l\in\cx^+_{\bJ}$ we denote by $\Pi(\l)$ the stabilizer of $[{}^\l\x]$ in $G$. Now
$\th$ is injective since if $\l\in\cx^+_J$, a subgroup $\Pi\in\cp^J$ is uniquely determined by the $\Pi$-stable
line in ${}^\l V$. Now let $(x_\l)\in\cc^*$. We show that the $H$-orbit of $(x_\l)$ is in $\th(\cp^J)$.
Let $\l\in\cx^+_{\bJ}$. We have $\G(x_{2\l})=E(x_\l,x_\l)$. Thus, $E_{x_\l,x_\l}$ is contained
in the irreducible summand of ${}^\l V\ot{}^\l V$ which is isomorphic to ${}^{2\l}V$, hence by a theorem of 
Kostant (see \cite{Ga} for the finite type case and \cite{PK} for the general case), we must have 
$[x_\l]=g_\l[{}^\l\x]$ for some 
$g_\l\in G$. Since $(x_\l)\in\cc^*$, for $\l,\l'$ in $\cx^+_{\bJ}$ we have 
$$\bE([g_{\l+\l'}({}^\l\x)],[g_{\l+\l'}({}^{\l'}\x)])=\bar\G([g_{\l+\l'}({}^{\l+\l'}\x)])=
\bE([g_\l({}^\l\x)],[g_{\l'}({}^{\l'}\x)]).$$
Since $\bE$ is injective, it follows that
$$[g_{\l+\l'}({}^\l\x)]=[g_\l({}^\l\x)],[g_{\l+\l'}({}^{\l'}\x)]=[g_{\l'}({}^{\l'}\x)],$$ 
so that 

(a) $g_\l\i g_{\l+\l'}\in\Pi(\l)$.
\nl
Assuming that $\l,\l'\in\cx^+_J$, we see that $g_\l\i g_{\l+\l'}\in\Pi^J$ and similarly 
$g_{\l'}\i g_{\l+\l'}\in\Pi^J$, so that $g_{\l'}\i g_\l\in\Pi^J$. Thus, there exists $g\in G$ such that for 
any $\l\in\cx^+_J$ we have
$g_\l=gp_\l$ with $p_\l\in\Pi^J$. Replacing $g_\l$ by $g_\l p_\l\i$, we see that we can assume that

(b) $g_\l=g$ for any $\l\in\cx^+_J$. 
\nl
If $\l\in\cx^+_{\bJ}$, $\l'\in\cx^+_J$, we have $\l+\l'\in\cx^+_J$ hence by (b),
$g_{\l+\l'}=g$, so that (a) implies $g_\l\i g\in\Pi(\l)$ and $[x_\l]=[g_\l({}^\l\x)]=[g({}^\l\x)]$.
Thus for any $\l\in\cx^+_{\bJ}$ we have  $x_\l=z_\l g({}^\l\x)$ for some $z_\l\in\CC^*$.
Since $(x_\l)\in\cc^*$ and $(g({}^\l\x))\in\cc^*$, we necessarily have $(z_\l)\in H$. Thus
the $H$-orbit of $(x_\l)$ is in the image of $\th$. The lemma is proved.

\subhead 1.3\endsubhead
Let $\cd$ be the category whose objects are pairs $(V,\b)$ where $V$ is a $\CC$-vector space and $\b$ is
a basis of $V$; a morphism from $(V,\b)$ to $(V',\b')$ is a $\CC$-linear map $f:V@>>>V'$ such that
for any $b\in\b$ we have $f(b)=\sum_{b'\in\b'}c_{b,b'}b'$ where $c_{b,b'}\in\NN$ for all $b,b'$ and
$c_{b,b'}=0$ for all but finitely many $b'$.

Let $K$ be a semifield. As in \cite{L19b} we define $K^!=K\sqc\{\circ\}$ where $\circ$ is a symbol.
We extend the sum and product on $K$ to a sum and product on $K^!$ by definining
$\circ+a=a$, $a+\circ=a$, $\circ\T a=\circ, a\T\circ=\circ$ for $a\in K$ and
$\circ+\circ=\circ, \circ\T\circ=\circ$. Thus $K^!$ becomes a monoid under addition and a monoid
under multiplication. Moreover the distributivity law holds in $K^!$.

A $K$-semivector space is an abelian (additive) semigroup $\cv$ with neutral element $\un\circ$
in which a map $K^!\T\cv@>>>\cv$, $(k,v)\m kv$ (``scalar multiplication'') is given 
such that $(kk')v=k(k'(v))$, $(k+k')v=kv+k'v$ for 
$k,k'$ in $K^!$, $v\in\cv$ and $k(v+v')=kv+kv'$ for $k\in K^!$, $v,v'$ in $\cv$; moreover we assume that
$k\un\circ=\un\circ$ for $k\in K^!$.

\mpb

Let $\cd(K)$ be the category whose objects are $K$-semivector spaces $\cv$; a morphism from $\cv$ to $\cv'$
is a map $f:\cv@>>>\cv'$ of semigroups preservng the neutral elements and commuting with scalar multiplication.
For any $\cv\in\cd(K)$ let $\End(\cv)=\Hom_{\cd(K)}(\cv,\cv)$; this is a monoid under composition of maps.

For $(V,\b)\in\cd$ let $V(K)$ be the set of formal sums $\x=\sum_{b\in\b}\x_bb$ with $\x_b\in K^!$ 
for all $b\in\b$  and $\x_b=\circ$ for all but finitely many $b$. We can define addition on $V(K)$
by $(\sum_{b\in\b}\x_bb)+(\sum_{b\in\b}\x'_bb)=\sum_{b\in\b}(\x_b+\x'_b)b$. We can define scalar 
multiplication by elements in $K^!$ by $k(\sum_{b\in\b}\x_b)=\sum_{b\in\b}(k\x_b)b$. Then $V(K)$
becomes an object of $\cd(K)$. The neutral element for addition is $\un\circ=\sum_{b\in\b}\circ b$.
Let $f$ be a morphism from $(V,\b)$ to $(V',\b')$ in $\cd$. 
For $b\in\b$ we have $f(b)=\sum_{b'\in\b'}c_{b,b'}b'$ where $c_{b,b'}\in\NN$. 
We define a map $f(K):V(K)@>>>V'(K)$ by
$f(K)(\sum_{b\in\b}\x_bb)=\sum_{b'\in\b'}(\sum_{b\in\b}c_{b,b'}\x_b)b'$.
Here for $c\in\NN,k\in K^!$ we set $ck=k+k+\do+k$ ($c$ terms) if $c>0$ and $ck=\circ$ if $c=0$.
Note that $f(K)$ is a morphism in $\cd(K)$. We have thus defined a functor $(V,\b)\m V(K)$ from $\cd$ to 
$\cd(K)$.

\mpb

Let $\l\in\cx^+$. We have $({}^\l V,{}^\l\b)\in\cd$ hence ${}^\l V(K)\in\cd(K)$ is defined. 
For $i\in I,m\in\ZZ$, the linear maps $e_i^{(n)},f_i^{(n)}$ from ${}^\l V$ to ${}^\l V$ (as in
\cite{L19b, 1.4}) are morphisms in $\cd$ (we use the positivity property \cite{L93, 22.1.7}
of ${}^\l\b$; in {\it loc.cit.} this property is stated assuming that the Cartan matrix is of
simply laced type, but the same proof applies in our case).
Hence they define morphisms $e_i^{(n)}(K),f_i^{(n)}(K)$ from ${}^\l V(K)$ to ${}^\l V(K)$.
For $i\in I,k\in K$ we define $i^k\in\End({}^\l V(K))$, $(-i)^k\in\End({}^\l V(K))$ by
$$i^k(b)=\sum_{n\in\NN}k^ne_i^{(n)}(K)b,\qua (-i)^k(b)=\sum_{n\in\NN}k^nf_i^{(n)}(K)b,$$ 
for any $b\in{}^\l\b$. 

For any $i\in I$ there is a well defined function $l_i:{}^\l\b@>>>\ZZ$ such that for $b\in{}^\l\b$, 
$t\in\CC^*$ we have $i(t)b=t^{l_i(b)}b$. (Here $i$ is viewed as a simple coroot homomorphism $\CC@>>>T$.)
For $i\in I,k\in K$ we define $\un i^k\in\End({}^\l V(K))$ by $\un i^k(b)=k^{l_i(b)}b$ for any $b\in{}^\l\b$.
As in \cite{L19b, 1.5}, the elements $i^k,(-i)^k,\un i^k$ (with $i\in I,k\in K$) in $\End({}^\l V(K))$
satisfy the relations in \cite{L19a, 2.10(i)-(vii)} defining the monoid $\fG(K)$ hence they define a monoid 
homomorphism $\fG(K)@>>>\End({}^\l V(K))$.
It follows that $\fG(K)$ acts on ${}^\l V(K)$.

\subhead 1.4\endsubhead
In the setup of 1.3 for $\l,\l'$ in $\cx^+$ we can view ${}^\l V\ot{}^{\l'}V$ with its basis 
$\cs={}^\l\b\ot{}^{\l'}\b$ as an object of $\cd$. Hence  $({}^\l V\ot{}^{\l'}V)(K)\in\cd(K)$ is defined.
We define $E(K):{}^\l V(K)\T{}^{\l'}V(K)@>>>({}^\l V\ot{}^{\l'}V)(K)$ by
$$(\sum_{b\in{}^\l\b}\x_bb),(\sum_{b'\in{}^{\l'}\b}\x'_{b'}b')\m\sum_{(b,b')\in\cs}\x_b\x'_{b'}(b\ot b').$$
(This is not a morphism in $\cd(K)$.) We define a map 
$$\End({}^\l V(K))\T\End({}^{\l'}V(K))@>>>\End(({}^\l V\ot{}^{\l'}V)(K))$$
by $(\t,\t')\m[b\ot b')\m E(K)(\t(b),\t'(b'))]$. Composing this map with the map
$$\fG(K)@>>>\End({}^\l V(K))\T\End({}^{\l'}V(K))$$
whose components are the maps 
$$\fG(K)@>>>\End({}^\l V(K)),\qua\fG(K)@>>>\End({}^{\l'}V(K))$$
 in 1.3 we obtain a map $\fG(K)@>>>\End(({}^\l V\ot{}^{\l'}V)(K))$ which 
is a monoid homomorphism. Thus $\fG(K)$ acts on $({}^\l V\ot{}^{\l'}V)(K)$; it also acts on
${}^\l V(K)\T{}^{\l'}V(K)$ (by 1.3) and the two actions are compatible with $E(K)$.

\mpb

Let $\G:{}^{\l+\l'}V@>>>{}^\l V\ot{}^{\l'}V$ be as in 1.1. For $b\in{}^{\l+\l'}\b$ we have 
$$\G(b)=\sum_{(b_1,b'_1)\in\cs}e_{b,b_1,b'_1}b_1\ot b'_1$$
where $e_{b,b_1,b'_1}\in\NN$. (This can be deduced from the positivity property \cite{L93, 14.4.13(b)} of the
homomorphism $r$ in \cite{L93, 1.2.12}.) Thus $\G$ is a morphism in $\cd$ hence 
$\G(K):{}^{\l+\l'}V(K)@>>>({}^\l V\ot{}^{\l'}V)(K)$ is a well defined morphism in $\cd(K)$.
Note that $\G(K)$ is compatible with the action of $\fG(K)$ on the two sides.

\subhead 1.5\endsubhead
In the setup of 1.3 let $\cc(K)$ be the set of all collections $\{x_\l\in{}^\l V(K);\l\in\cx^+_{\bJ}\}$
such that for any $\l,\l'$ in $\cx^+_{\bJ}$ we have $\G(K)(x_{\l+\l'})=E(K)(x_\l,x_{\l'})$.
Let $\cc^*(K)$ be the set of all $(x_\l)\in\cc(K)$ such that $x_\l\ne\un\circ$ for any $\l\in\cx^+_{\bJ}$.
Let $H(K)$ be the group (multiplication component by component)
consisting of all collections $\{z_\l\in K;\l\in\cx^+_{\bJ}\}$
such that for any $\l,\l'$ in $\cx^+_{\bJ}$ we have $z_{\l+\l'}=z_\l z_{\l'}$.
Now $H(K)$ acts on $\cc(K)$ by $(z_\l),(x_\l)\m(z_\l x_\l)$.
This restricts to a free action of $H(K)$ on $\cc^*(K)$. Let $\cp^J(K)$ be the set of orbits for this action.
Note that $\fG(K)$ acts on $\cc(K)$ by acting component by component (see 1.3); we use that $E(K),\G(K)$ are 
compatible with the $\fG(K)$-actions (see 1.4). This induces a $\fG(K)$-action on $\cp^J(K)$.

\subhead 1.6\endsubhead
In this subsection we assume that $K=\RR_{>0}$. If $(x_\l)\in\cc^*(K)$, we can view $(x_\l)$ as an element 
of $\cc^*$ by viewing ${}^\l V(K)$ as a subset of ${}^\l V$ in an obvious way.
The inclusion $\cc^*(K)\sub\cc^*$ is compatible with the actions of $H(K)$ and $H$
(we have $H(K)\sub H$) hence it induces an (injective) map $\cp^J(K)@>>>{}'\cp^J$. Composing this with the 
inverse of the bijection $\cp^J@>>>{}'\cp^J$ (see 1.2) we obtain an injective map $\cp^J(K)@>>>\cp^J$.
We define $\cp^J_{\ge0}$ to be the image of this map.

Assuming further that our Cartan matrix is of finite type, we show that the last definition of $\cp^J_{\ge0}$
agrees with the definition in \cite{L98}. Applying 
\cite{L98, 3.4} to a $\l\in\cx^+_J$ with large enough coordinates we see that $\cp^J_{\ge0}$ (in the new
definition) is contained in $\cp^J_{\ge0}$ (in the definition of \cite{L98}). 
The reverse inclusion follows from \cite{L98, 3.2}. 

\subhead 1.7\endsubhead
Any homomorphism of semifields $K@>>>K'$ induces in an obvious way a map $\cp^J(K)@>>>\cp^J(K')$.

\subhead 1.8\endsubhead
We expect that when $K'$ is the semifield $\{1\}$ with one element, 
one can identify $\cp^\emp(K')$ with the set of 
pairs $(a,a')$ in the Weyl group $W$ of $G$ such that $a\le a'$ for the standard partial order of $W$.
If $K$ is any semifield one can also expect that the fibre of the map $\cp^\emp(K)@>>>\cp^\emp(\{1\})$
induced by the obvious map $K@>>>\{1\}$ (see 1.7) at the element corresponding to $(a,a')$ is in bijection 
with $K^{|a'|-|a|}$ where $a\m|a|$ is the length function on $W$.

\head 2. The semiring $M(K)$\endhead
\subhead 2.1\endsubhead
In this section we assume that our Cartan matrix is of finite type.
Let $K$ be a semifield. Let $M(K)=\op_{\l\in\cx^+_{\bJ}}{}^\l V(K)$ viewed as a monoid 
under addition and with scalar multiplication by elements of $K^!$.

We define a multiplication $\mu:M(K)\T M(K)@>>>M(K)$ which is ``bilinear'' with respect to
addition and scalar multiplication and satisfies $\mu(b_1,b'_1)=\sum_{b\in{}^{\l+\l'}\b}e_{b,b_1,b'_1}b$
where $\l\in\cx^+_{\bJ},\l'\in\cx^+_{\bJ}$, $b_1\in{}^\l\b,b'_1\in{}^{\l'}\b$ and $e_{b,b_1,b'_1}\in \NN$
(viewed as an element of $K^!$) is as in the definition of $\G(K)$ in 1.4, so that it comes from the 
homomorphism $r$ in \cite{L93, 1.2.12}. This can be viewed as a direct sum of ``transposes'' of maps like 
$\G(K)$.
From the properties of $r$ we see that the multiplication $\mu$ is associative and commutative; it is clearly
distributive with respect to addition. This multiplication has a unit element, given by the unique element in 
$\b^\l$ with $\l=0$. Note that $M(K)$ is a semiring.
Now 
$M(K)$ can be viewed as a form over $K$ of the coordinate ring of $G/U^+$ where $U^+$ is the unipotent radical
of $B^+$.
Let $M'(K)$ be the set of maps $M(K)@>>>K^!$ which are compatible with addition, multiplication and with
scalar multiplication by elements of $K^!$, take the unit element of $M(K)$ to the unit element of $K^!$ and take
the element with all components equal to $\un\circ$ to $\circ\in K^!$.
It is easy to show that $M'(K)$ is in canonical bijection with $\cc(K)$.

\enddocument